\documentclass{article}%
\usepackage{amssymb}
\usepackage{amsfonts}
\usepackage{amsmath}
\usepackage{graphicx}%
\providecommand{\U}[1]{\protect\rule{.1in}{.1in}}

\begin{document}

\title{Forward-backward systems of stochastic differential equations generated by
Bernstein diffusions }
\author{Ana Bela Cruzeiro$^{\ast}$ and Pierre-A. Vuillermot$^{\ast\ast}$\\Grupo de F\'{\i}sica Matem\'{a}tica da Universidade de Lisboa$^{\ast}$\\UMR-CNRS 7502, Institut \'{E}lie Cartan de Lorraine, France$^{\ast\ast}$}
\date{}
\maketitle

\begin{abstract}
In this short article we present new results that bring about hitherto unknown
relations between certain Bernstein diffusions wandering in bounded convex
domains of Euclidean space on the one hand, and processes which typically
occur in forward-backward systems of stochastic differential equations on the
other hand. A key point in establishing such relations is the fact that the
Bernstein diffusions we consider are actually reversible It\^{o} diffusions.

\end{abstract}

\section{Statement of the main results}

The theory of Bernstein processes goes back to \cite{bernstein} which
elaborates on the seminal contribution that was set forth in the very last
section of \cite{schroedinger}. It was subsequently thoroughly developed in
\cite{jamison}, and has ever since played an important r\^{o}le in probability
theory and in various areas of mathematical physics as testified for instance
by the works \cite{carlen}-\cite{cruzeirozambrini}, \cite{gulicasteren},
\cite{prizamb}, \cite{vuizamb}, \cite{zambrini} and the many references
therein. There are several equivalent ways to define a Bernstein process, but
we settle here for a variant which is tailored to our needs. Let
$D\subset\mathbb{R}^{d}$ be a bounded open convex subset whose smooth boundary
$\partial D$ is $\mathcal{C}^{2+\alpha}$ with $\alpha\in\left(  0,1\right)  $;
let us write $\overline{D}:=D\cup\partial D$ and let $T\in\left(
0,+\infty\right)  $ be arbitrary.

\bigskip

\textbf{Definition 1.1.} We say the $\overline{D}$-valued process $Z_{\tau
\in\left[  0,T\right]  }$ defined on the complete probability space $\left(
\Omega,\mathcal{F},\mathbb{P}\right)  $ is a Bernstein process if the
following conditional expectations satisfy the relation%
\begin{equation}
\mathbb{E}\left(  h(Z_{r})\left\vert \mathcal{F}_{s}^{+}\vee\mathcal{F}%
_{t}^{-})\right.  \right)  =\mathbb{E}\left(  h(Z_{r})\left\vert Z_{s}%
,Z_{t}\right.  \right)  \label{conditionalexpectations}%
\end{equation}
for every Borel mesurable function $h:\overline{D}\rightarrow\mathbb{R}$, and
for all $r,s,t$ satisfying $r\in(s,t)\subset\left[  0,T\right]  $. In
(\ref{conditionalexpectations}), $\mathcal{F}_{s}^{+}$ denotes the $\sigma
$-algebra generated by the $Z_{\tau}$'s for all $\tau\in\left[  0,s\right]  $,
while $\mathcal{F}_{t}^{-}$ is that generated by the $Z_{\tau}$'s for all
$\tau\in\left[  t,T\right]  $.

\bigskip

According to \cite{jamison} Bernstein processes are neither Markovian nor
reversible in general, but there are special classes of them which are
reversible It\^{o} diffusions as we shall see below. Since this conveys the
idea that such diffusions can evolve forward and backward in time while
satisfying It\^{o} equations for some suitably constructed Wiener processes,
they constitute the bare minimum we need to establish relations with certain
forward-backward systems of stochastic differential equations. Following
\cite{bensoussan}, \cite{bismut} and \cite{parpeng} which were written in the
context of stochastic optimal control, the theory of such systems has been
considerably developed in recent times and has led to various applications in
mathematical finance and partial differential equations (see e.g.
\cite{delarue}, \cite{dutang}, \cite{macvitanic}, \cite{majong},
\cite{parzhang} and their references). In this Note we adopt the standard
definition for them which can be found for instance in \cite{delarue} or
\cite{majong}, up to the essential difference that the Wiener processes which
occur in our context are not \textit{a priori} given but rather determined
\textit{a posteriori }by\textit{\ }the Bernstein diffusions themselves. This
leads to the following notion, where%
\begin{equation}
f,g:\overline{D}\times\mathbb{R}^{d}\times\left[  0,T\right]  \mapsto
\mathbb{R}^{d} \label{specialvectorfields}%
\end{equation}
are continuous functions, and where we write $X=\left(  X_{1},...,X_{d}%
\right)  $ for any vector $X\in\mathbb{R}^{d^{2}}$with each $X_{i}%
\in\mathbb{R}^{d}$ and $\left\vert .\right\vert $ for the Euclidean norm in
$\mathbb{R}^{d}$ or $\mathbb{R}^{d^{2}}$:

\bigskip

\textbf{Definition 1.2.} We say the $\overline{D}\times\mathbb{R}^{d}%
\times\mathbb{R}^{d^{2}}$-valued process $\left(  A_{\tau},B_{\tau},C_{\tau
}\right)  _{\tau\in\left[  0,T\right]  }$ defined on $\left(  \Omega
,\mathcal{F},\mathbb{P}\right)  $ is a solution to a forward-backward system
of stochastic differential equations with initial condition $\eta$ and final
condition $\kappa$ if the following four conditions hold:

(a) The process $W_{\tau\in\left[  0,T\right]  }$ defined by%
\[
W_{\tau}:=A_{\tau}-A_{0}-\int_{0}^{\tau}d\sigma f\left(  A_{\sigma},B_{\sigma
},\sigma\right)
\]
is a Wiener process on $\left(  \Omega,\mathcal{F},\mathbb{P}\right)  $
relative to its natural increasing filtration $\mathcal{F}_{\tau\in\left[
0,T\right]  }$, and $\left(  A_{\tau},B_{\tau},C_{\tau}\right)  _{\tau
\in\left[  0,T\right]  }$ is progressively measurable with respect to
$\mathcal{F}_{\tau\in\left[  0,T\right]  }$.

(b) We have%
\begin{equation}
A_{t}=\eta+\int_{0}^{t}d\tau f\left(  A_{\tau},B_{\tau},\tau\right)  +W_{t}
\label{abstractforward1}%
\end{equation}
$\mathbb{P}$-a.s. for every $t\in\left[  0,T\right]  $, where $\eta
:\Omega\mapsto D$ is $\mathcal{F}_{0}$-measurable and satisfies $\mathbb{E}%
\left\vert \eta\right\vert ^{2}<+\infty$.

(c) We have%
\begin{equation}
B_{t,i}=\kappa_{i}-\int_{t}^{T}d\tau g_{i}\left(  A_{\tau},B_{\tau}%
,\tau\right)  -\int_{t}^{T}\left(  C_{\tau,i},dW_{\tau}\right)  _{\mathbb{R}%
^{d}} \label{abstractbackward2}%
\end{equation}
$\mathbb{P}$-a.s. for every $t\in\left[  0,T\right]  $ and every $i\in\left\{
1,...,d\right\}  $, where $\kappa:\Omega\mapsto\mathbb{R}^{d}$ is
$\mathcal{F}_{T}$-measurable and satisfies $\mathbb{E}\left\vert
\kappa\right\vert ^{2}<+\infty$. In (\ref{abstractbackward2}),
$(.,.)_{\mathbb{R}^{d}}$ stands for the Euclidean inner product in
$\mathbb{R}^{d}$ and the second integral is the forward It\^{o} integral
defined with respect to $\mathcal{F}_{\tau\in\left[  0,T\right]  }$, which we
assume to be well-defined.

(d) We have%
\[
\mathbb{E}\int_{0}^{T}d\tau\left(  \left\vert A_{\tau}\right\vert
^{2}+\left\vert B_{\tau}\right\vert ^{2}+\left\vert C_{\tau}\right\vert
^{2}\right)  <+\infty.
\]

\bigskip

In order to show that there exist special classes of Bernstein processes which
generate forward-backward systems in the above sense, we now consider
parabolic initial-boundary value problems\textit{\ }of the form%
\begin{align}
\partial_{t}u(x,t)  &  =\frac{1}{2}\bigtriangleup_{x}u(x,t)-\left(
l(x,t),\nabla_{x}u(x,t)\right)  _{\mathbb{R}^{d}}-V(x,t)u(x,t),\nonumber\\
\text{\ }(x,t)  &  \in D\times\left(  0,T\right]  ,\nonumber\\
u(x,0)  &  =\varphi(x),\text{ \ \ }x\in D,\nonumber\\
\frac{\partial u(x,t)}{\partial n(x)}  &  =0,\text{ \ \ \ }(x,t)\in\partial
D\times\left(  0,T\right]  \label{parabolicproblem}%
\end{align}
along with the corresponding adjoint final-boundary value problems%
\begin{align}
-\partial_{t}v(x,t)  &  =\frac{1}{2}\bigtriangleup_{x}%
v(x,t)+\operatorname{div}_{x}\left(  v(x,t)l(x,t)\right)
-V(x,t)v(x,t),\nonumber\\
(x,t)  &  \in D\times\left[  0,T\right)  ,\nonumber\\
v(x,T)  &  =\psi(x),\text{ \ \ }x\in D,\nonumber\\
\frac{\partial v(x,t)}{\partial n(x)}  &  =0,\text{ \ \ \ }(x,t)\in\partial
D\times\left[  0,T)\right.  . \label{adjointproblem}%
\end{align}
In the preceding relations $n(x)$ denotes the unit outer normal vector at
$x\in\partial D$, $l$ is an $\mathbb{R}^{d}$-valued vector-field and $V$,
$\varphi$, $\psi$ are real-valued functions which satisfy the following
hypotheses, respectively:

\bigskip

(L) For $l:\overline{D}\times\left[  0,T\right]  \mapsto\mathbb{R}^{d}$ we
have $l_{i}\in\mathcal{C}^{2+\alpha,\frac{\alpha}{2}}(\overline{D}%
\times\left[  0,T\right]  )$ for every $i\in\left\{  1,...,d\right\}  $.

(V) The function $V:\overline{D}\times\left[  0,T\right]  \mapsto\mathbb{R}$
is such that $V\in\mathcal{C}^{1,\frac{\alpha}{2}}(\overline{D}\times\left[
0,T\right]  )$.

(IF) We have $\varphi,\psi\in\mathcal{C}^{2+\alpha}(\overline{D})$ with
$\varphi>0$ and $\psi>0$ satisfying Neumann's boundary condition, that is,%
\[
\frac{\partial\varphi(x)}{\partial n(x)}=\frac{\partial\psi(x)}{\partial
n(x)}=0,\text{ \ \ \ }x\in\partial D.
\]

\bigskip

Under such conditions there exists a Bernstein process $Z_{\tau\in\left[
0,T\right]  }$ wandering in $\overline{D}$ which turns out to be a reversible
It\^{o} diffusion according to the theory developed in \cite{vuizamb}, to
which we refer the reader for details. What this means is that $Z_{\tau
\in\left[  0,T\right]  }$ may be considered simultaneously either as a forward
or as a backward Markov diffusion on some probability space $\left(
\Omega,\mathcal{F},\mathbb{P}_{\mu}\right)  $, where $\mathbb{P}_{\mu}$ is
uniquely determined by the positive function%
\begin{equation}
\mu(E\times F):=\int_{E\times F}dxdy\varphi(x)g(y,T;x,0)\psi(y)
\label{jointmeasure}%
\end{equation}
defined for all $E,F\in\mathcal{B}\left(  \overline{D}\right)  $ with
$\mathcal{B}\left(  \overline{D}\right)  $ the Borel $\sigma$-algebra on
$\overline{D}$. In the preceding expression $g$ stands for the parabolic Green
function associated with (\ref{parabolicproblem}), which is well-defined and
positive for all $x,y\in\overline{D}$. In addition, (\ref{jointmeasure}) must
satisfy%
\[
\int_{D\times D}dxdy\varphi(x)g(y,T;x,0)\psi(y)=1.
\]
Furthermore the drift of $Z_{\tau\in\left[  0,T\right]  }$ is%
\[
b^{\ast}(x,t)=l(x,t)+\nabla_{x}\ln v_{\psi}(x,t)
\]
in the forward case, and%
\[
b(x,t)=l(x,t)-\nabla_{x}\ln u_{\varphi}(x,t)
\]
in the backward case, where $u_{\varphi}$ is the unique positive classical
solution to (\ref{parabolicproblem}) and $v_{\psi}$ the unique classical
positive solution to (\ref{adjointproblem}). Finally $Z_{\tau\in\left[
0,T\right]  }$ satisfies two It\^{o} equations in the weak sense, namely, the
forward equation%
\begin{equation}
Z_{t}=Z_{0}+\int_{0}^{t}d\tau b^{\ast}\left(  Z_{\tau},\tau\right)
+W_{t}^{\ast} \label{itoequation1}%
\end{equation}
and the backward equation%
\begin{equation}
Z_{t}=Z_{T}-\int_{t}^{T}d\tau b\left(  Z_{\tau},\tau\right)  +W_{t}
\label{itoequation2}%
\end{equation}
$\mathbb{P}_{\mu}$-a.s. for every\textit{\ }$t\in\left[  0,T\right]  $, for
two suitably constructed $d$-dimensional Wiener processes $W_{\tau\in\left[
0,T\right]  }^{\ast}$ and $W_{\tau\in\left[  0,T\right]  }$.

It is precisely all these properties along with some more refinements that
lead to the desired relations with forward-backward systems. For instance, let
us choose the functions in (\ref{specialvectorfields}) as%
\begin{equation}
f(x,y,t)=l\left(  x,t\right)  +y \label{specialvectorfield1}%
\end{equation}
and%
\begin{equation}
g_{i}(x,y,t)=\partial_{x_{i}}(V\left(  x,t\right)  -\operatorname{div}%
_{x}l(x,t))-\left(  y,\nabla_{x}l_{i}(x,t)\right)  _{\mathbb{R}^{d}}
\label{specialvectorfield2}%
\end{equation}
for every\textit{\ }$i\in\left\{  1,...,d\right\}  $. One of our typical
results is then the following:

\bigskip

\textbf{Theorem.} \textit{Let us assume that Hypotheses (L), (V) and (IF)
hold, and that the vector-field }$l$ \textit{is conservative in }$D$.
\textit{Then, there exist} \textit{a probability measure} $\mu$ on
$\mathcal{B}\left(  \overline{D}\right)  \times\mathcal{B}\left(  \overline
{D}\right)  $, \textit{a probability space }$\left(  \Omega,\mathcal{F}%
,\mathbb{P}_{\mu}\right)  $\textit{\ and a} $\overline{D}\times\mathbb{R}%
^{d}\times\mathbb{R}^{d^{2}}$\textit{-valued} \textit{process} $\left(
A_{\tau},B_{\tau},C_{\tau}\right)  _{\tau\in\left[  0,T\right]  }$ \textit{on}
$\left(  \Omega,\mathcal{F},\mathbb{P}_{\mu}\right)  $,\textit{\ such that the
following properties hold:}

\textit{(a) The process }$\left(  A_{\tau},B_{\tau},C_{\tau}\right)  _{\tau
\in\left[  0,T\right]  }$ \textit{is the unique solution to the
forward-backward system whose stochastic differential equations are of the
form (\ref{abstractforward1}) and (\ref{abstractbackward2}) with }$f$
\textit{and }$g$ \textit{given by (\ref{specialvectorfield1}) and
(\ref{specialvectorfield2}), respectively. The uniqueness is meant as
uniqueness in law. }

\textit{(b) The initial random vector }$\eta:\Omega\mapsto D$\textit{\ has the
distribution density }$x\mapsto\varphi(x)v_{\psi}(x,0)$ \textit{for all} $x\in
D$, \textit{while }$\kappa:\Omega\mapsto\mathbb{R}^{d}$ \textit{is of the
form}%
\[
\kappa=\nabla_{x}\ln\psi\left(  A_{T}\right)  \text{.}%
\]

\textit{(c) The process }$A_{\tau\in\left[  0,T\right]  }$ \textit{is a
reversible Bernstein diffusion in} $\overline{D}$.

\section{A brief sketch of the proof}

We first consider the translated forward drift%
\begin{equation}
c^{\ast}(x,t):=b^{\ast}(x,t)-l(x,t)=\nabla_{x}\ln v_{\psi}(x,t)
\label{translatedforwardrift}%
\end{equation}
for all $(x,t)\in\overline{D}\times\left[  0,T\right]  $, and then rewrite
(\ref{itoequation1}) as
\begin{equation}
Z_{t}=Z_{0}+\int_{0}^{t}d\tau\left(  c^{\ast}\left(  Z_{\tau},\tau\right)
+l\left(  Z_{\tau},\tau\right)  \right)  +W_{t}^{\ast}.
\label{newitoequation1}%
\end{equation}
We can also verify that the process $c^{\ast}\left(  Z_{\tau},\tau\right)
_{\tau\in\left[  0,T\right]  \text{ \ }}$satisfies the backward It\^{o}
equation%
\begin{align}
&  c_{i}^{\ast}\left(  Z_{t},t\right) \nonumber\\
&  =c_{i}^{\ast}\left(  Z_{T},T\right)  -\int_{t}^{T}d\tau\left(
\partial_{x_{i}}(V\left(  Z_{\tau},\tau\right)  -\operatorname{div}%
_{x}l(Z_{\tau},\tau))-\left(  c^{\ast}(Z_{\tau},\tau),\nabla_{x}l_{i}(Z_{\tau
},\tau)\right)  _{\mathbb{R}^{d}}\right) \nonumber\\
&  -\int_{t}^{T}\left(  \nabla_{x}c_{i}^{\ast}\left(  Z_{\tau},\tau\right)
,d^{+}W_{\tau}^{\ast}\right)  _{\mathbb{R}^{d}} \label{backwarditoequation}%
\end{align}
$\mathbb{P}_{\mu}$-a.s. for every $t\in\left[  0,T\right]  $ and every
$i\in\left\{  1,...,d\right\}  $, where the second integral in
(\ref{backwarditoequation}) is the forward It\^{o} integral defined with
respect to $W_{\tau\in\left[  0,T\right]  }^{\ast}$ and its natural increasing
filtration $\mathcal{F}_{\tau\in\left[  0,T\right]  }^{+}$. It is then
sufficient to choose $\left(  A_{\tau},B_{\tau},C_{\tau}\right)  _{\tau
\in\left[  0,T\right]  }=\left(  Z_{\tau},c^{\ast}\left(  Z_{\tau}%
,\tau\right)  ,\nabla c^{\ast}\left(  Z_{\tau},\tau\right)  \right)  _{\tau
\in\left[  0,T\right]  }$, where we have written%
\begin{equation}
\nabla c^{\ast}\left(  Z_{\tau},\tau\right)  _{\tau\in\left[  0,T\right]
}:=\left(  \nabla_{x}c_{1}^{\ast}\left(  Z_{\tau},\tau\right)  ,...,\nabla
_{x}c_{d}^{\ast}\left(  Z_{\tau},\tau\right)  \right)  _{_{\tau\in\left[
0,T\right]  }}. \label{gradientprocess}%
\end{equation}
This choice shows indeed that we can identify (\ref{abstractforward1}) and
(\ref{abstractbackward2}) with (\ref{newitoequation1}) and
(\ref{backwarditoequation}), respectively, if we take
(\ref{specialvectorfield1}) and (\ref{specialvectorfield2}) into account. The
remaining statements follow from simple considerations.

We can obtain so to speak a dual result when we consider the process $\left(
A_{\tau},B_{\tau},C_{\tau}\right)  _{\tau\in\left[  0,T\right]  }=\left(
Z_{\tau},c\left(  Z_{\tau},\tau\right)  ,\nabla c\left(  Z_{\tau},\tau\right)
\right)  _{\tau\in\left[  0,T\right]  }$, where
\begin{equation}
c(x,t):=b(x,t)-l(x,t)=-\nabla_{x}\ln u_{\varphi}(x,t)
\label{translatedbackwardrift}%
\end{equation}
is the translated backward drift associated with $Z_{\tau\in\left[
0,T\right]  }$, which we consider this time as a backward It\^{o} diffusion
satisfying (\ref{itoequation2}) or, equivalently,%
\begin{equation}
Z_{t}=Z_{T}-\int_{t}^{T}d\tau\left(  c(Z_{\tau},\tau)+l(Z_{\tau},\tau)\right)
+W_{t}. \label{newitoequation2}%
\end{equation}
Indeed, in this case it is the process $c\left(  Z_{\tau},\tau\right)
_{\tau\in\left[  0,T\right]  \text{ }}$ that satisfies a forward It\^{o}
equation. We refer the reader to the forthcoming typescript \cite{cruvui} for
details and complete proofs.

\bigskip

\bigskip

\bigskip

\bigskip

\bigskip

\bigskip

\bigskip

\bigskip

\bigskip

\bigskip


\begin{thebibliography}{99}                                                                                               %


\bibitem {bensoussan}\textsc{A. Bensoussan, }\textit{Stochastic maximum
principle for distributed parameter systems,} Journal of the Franklin
Institute \textbf{315 }(1983) 387-406.

\bibitem {bernstein}\textsc{S. Bernstein, }\textit{Sur les liaisons entre les
grandeurs al\'{e}atoires,} in: Verhandlungen des Internationalen
Mathematikerkongress, \textbf{1} (1932) 288-309.

\bibitem {bismut}\textsc{J. M. Bismut,}\textit{\ Th\'{e}orie probabiliste du
contr\^{o}le des diffusions,} Memoirs of the AMS \textbf{176}, American
Mathematical Society, Providence, 1973.

\bibitem {carlen}\textsc{E. A. Carlen,} \textit{Conservative diffusions,
}Communications in Mathematical Physics \textbf{94} (1984) 293-315.

\bibitem {cruvui}\textsc{A. B. Cruzeiro, P.-A. Vuillermot,} \textit{Bernstein
diffusions and forward-backward stochastic differential equations, }(2013).

\bibitem {cruwuzamb}\textsc{A. B. Cruzeiro, L. Wu, J. C. Zambrini,
}\textit{Bernstein processes associated with Markov processes,} In: Rebolledo,
R. (ed.) Stochastic Analysis and Mathematical Physics, Birkh\"{a}user, Basel (2000).

\bibitem {cruzeirozambrini}\textsc{A. B. Cruzeiro, J. C. Zambrini,
}\textit{Malliavin calculus and Euclidean quantum mechanics, I. Functional
Calculus,} Journal of Functional Analysis \textbf{96} (1991) 62-95.

\bibitem {delarue}\textsc{F. Delarue, }\textit{On the existence and uniqueness
of solutions to FBSDEs in a non-degenerate case, }Stochastic Processes and
their Applications \textbf{99} (2002) 209-286.

\bibitem {dutang}\textsc{K. Du, S. Tang,} \textit{Strong solution of backward
stochastic partial differential equations in }$\mathcal{C}^{2}$%
\textit{-domains,} Probability Theory and Related Fields \textbf{154} (2012) 255-285.

\bibitem {gulicasteren}\textsc{A. Gulisashvili, J. A. van Casteren,}
\textit{Non-Autonomous Kato Classes and Feynman-Kac Propagators,} World
Scientific, Singapore (2000).

\bibitem {jamison}\textsc{B. Jamison, }\textit{Reciprocal processes,}
Zeitschrift f\"{u}r Wahrscheinlichkeitstheorie und Verwandte Gebiete
\textbf{30} (1974) 65-86.

\bibitem {macvitanic}\textsc{J. Ma, J. Cvitani\'{c},} \textit{Reflected
forward-backward SDEs and obstacle problems with boundary conditions, }Journal
of Applied Mathematics and Stochastic Analysis \textbf{14:2 }(2001) 113-138.

\bibitem {majong}\textsc{J. Ma, J. Yong,} \textit{Forward-Backward Stochastic
Differential Equations and their Applications,} Lectures Notes in Mathematics
\textbf{1702}, Springer, New York (2007).

\bibitem {parpeng}\textsc{E. Pardoux, S. G. Peng, }\textit{Adapted solution of
a backward stochastic differential equation,} Systems and Control Letters
\textbf{14} (1990) 55-61.

\bibitem {parzhang}\textsc{E. Pardoux, S. Zhang,} \textit{Generalized BSDEs
and nonlinear Neumann boundary-value problems,} Probability Theory and Related
Fields \textbf{110} (1998)\ 535-558.

\bibitem {prizamb}\textsc{N. Privault, J. C. Zambrini, }\textit{Markovian
bridges and reversible diffusion processes with jumps,} Annales de l'Institut
Henri Poincar\'{e} PR \textbf{40} (2004) 599-633.

\bibitem {schroedinger}\textsc{E. Schr\"{o}dinger, }\textit{Sur la th\'{e}orie
relativiste de l'\'{e}lectron et l'interpr\'{e}tation de la m\'{e}canique
quantique,} Annales de l'Institut Henri-Poincar\'{e} \textbf{2} (1932) 269-310.

\bibitem {vuizamb}\textsc{P.-A. Vuillermot, J. C. Zambrini,} \textit{Bernstein
diffusions for a class of linear parabolic partial differential equations,}
Journal of Theoretical Probability \textbf{DOI 10.1007/s10959-012-0426-3} (2012).

\bibitem {zambrini}\textsc{J. C. Zambrini,} \textit{Variational processes and
stochastic versions of mechanics,} Journal of Mathematical Physics \textbf{27}
(1986) 2307-2330.
\end{thebibliography}
\end{document}